\documentclass{article}

\usepackage{amsmath}   
\usepackage{amssymb}   
\usepackage{accents}
\usepackage{graphicx}
\usepackage{color}
\usepackage{tikz}
\usepackage{paralist}
\usepackage{enumerate} 
\usepackage{amsthm}  
\usepackage[english]{babel}

\theoremstyle{theorem}            
\newtheorem{theorem}{Theorem}[section]

\theoremstyle{definition}           

\newtheorem{proposition}[theorem]{Proposition}
\newtheorem{lemma}[theorem]{Lemma}
\newtheorem{corollary}[theorem]{Corollary}

\newtheorem{question}{Question}

\newcommand{\N}{\mathbb{N}}
\newcommand{\R}{\mathbb{R}}
\newcommand{\Z}{\mathbb{Z}}
\newcommand{\Q}{\mathbb{Q}}

\newcommand{\cantor}{\ensuremath{2^\N}}
\newcommand{\conj}{~\&~}
\newcommand{\ignore}[1]{}
\renewcommand{\subset}{\subseteq}
\newcommand{\ciso}{\cong_{\mathcal{L}}^c}
\newcommand{\iffdef}{\ \mathrel{\mathord{:\!\!\iff}}}
\newcommand{\bigdoublevee}{%
  \mathop{
    \mathchoice{\bigvee\mkern-15mu\bigvee}
               {\bigvee\mkern-12.5mu\bigvee}
               {\bigvee\mkern-12.5mu\bigvee}
               {\bigvee\mkern-11mu\bigvee}
    }
}

\title{Graphings of arithmetical equivalence relations}
\author{Tyler Arant\footnote{University of California, Los Angeles}}
\date{\today}

\begin{document}

\maketitle

\begin{abstract}
This paper studies when an arithmetical equivalence
relation $E$ can be realized as the connectedness
relation of a graph $G$ which is simpler to define
than $E$.  Several examples of such equivalence relations
are established.  In particular, it is proved that 
the $\Sigma^0_3$
relation of computable isomorphism of structures on $\N$
in a computable first-order language is $\Pi^0_2$-graphable, i.e.,
is the connectedness relation of a $\Pi^0_2$ graph.  
Graphings of Friedman-Stanley jumps are studied, including 
an arithmetical construction of a graphing
of the Friedman-Stanley jump of $E$ from a graphing of $E$.
\end{abstract}

\section{Introduction}

Let $\Gamma$ be a pointclass.\footnote{
A \textbf{pointclass} $\Gamma$ is an operation which assigns
every Polish (separable, completely metrizable)
space $X$ to a collection $\Gamma(X)$ of 
subsets of $X$.
}
An equivalence relation $E$ on a Polish space $X$ is called
\textbf{$\Gamma$-graphable} if there is a (simple, undirected)
graph $G$ in pointclass $\Gamma$ such that the connectedness equivalence
relation of $G$ is equal to $E$, i.e.,
\[
xEy \iff \text{there is a path in $G$ connecting $x$ to $y$}.
\]
We also say that such a $G$ is a \textbf{$\Gamma$ graphing}
of $E$.  

For $k\in \N$, a graph $G$ has diameter $k$ if every pair of
$G$-connected points $x, y$ are connected by a path of
length at most $k$ and, moreover, $k$ is the least integer with
this property.  We say $E$ is \textbf{$\Gamma$-graphable with diameter $k$}
if it has a $\Gamma$ graphing $G$ whose diameter is $k$.

The topic of $\Gamma$-graphability was initially studied 
in \cite{thesis} and greatly expanded upon in \cite{akl}.
Both of those works mainly studied which analytic equivalence
relations are Borel graphable.  
In this paper, we are interested in
arithmetical equivalence relations which have graphings with
simpler arithmetical definitions.\footnote{
By \textbf{arithmetical}, we mean belonging to one
of the (lightface) pointclasses $\Sigma^0_n$, $n\geq 1$. 
I.e., an arithmetical set is definable using a finite number
of quantifiers over $\N$ in front of a computable predicate.
}

Arithmetical equivalence
relations play an important role
in both descriptive set theory and computability theory.
We highlight some of the most important examples, with a focus
on those for which we will prove graphability results.

(A) The eventual equality 
equivalence relation $E_0$ on the Cantor space
$\cantor$ is 
defined by
\[
xE_0 y \iffdef (\exists m)(\forall n\geq m)[x(n)=y(n)].
\]
$E_0$ is clearly $\Sigma^0_2$, and it serves as a vital ``benchmark''
in the study of Borel equivalence relations.  The Glimm-Effros dichotomy
theorem from \cite{hkl} shows that continuously
embedding $E_0$ is the canonical obstruction to a
Borel equivalence relation being smooth.  Borel reducibility
to $E_0$ also provides
an alternative characterization of hyperfiniteness for
countable Borel equivalence relations, see \cite{djk}.

$E_0$ is the connectedness relation of an important $\Delta^0_2$ graph, which
is denoted $G_0$.  The $G_0$ dichotomy theorem (see \cite{kst}) establishes
that $G_0$ is the minimal obstruction to an analytic graph having
a countable Borel chromatic number.  The importance of both $E_0$ and
its graphing $G_0$ serves as motivation to understand the graphings
of arithmetical equivalence relations.  

(B) Turing equivalence, denoted $\equiv_T$, as an equivalence relation on 
$\cantor$ is a properly $\Sigma^0_3$ equivalence relation; 
in fact,
it is a (boldface) $\boldsymbol{\Sigma}^0_3$-complete 
(see \cite{rss}, Corollary 22).
We will see in Section \ref{basic} that $\equiv_T$ is 
$\Pi^0_2$-graphable with diameter $2$.  

(C) An \textbf{$m$-reduction} from $A\subset \N$ to  $B\subset \N$ is a 
computable total $f:\N\rightarrow \N$ such that 
$n\in A$ if and only if $f(n)\in B$ for all $n\in \N$.
We write $A\leq_m B$ when $A$ can be $m$-reduced to $B$.
\textbf{$m$-equivalence} is defined by
$A\equiv_m B$
if and only if $A\leq_m B$ and $B\leq_m A$.  

If an $m$-reduction
$f$ is a bijection, then it is called a 
\textbf{$1$-equivalence}.
$A$ and $B$ are \textbf{$1$-equivalent}, denoted $A\equiv_1 B$,
if there is a $1$-equivalence between them.\footnote{
Note that by a theorem of Myhill, $A\equiv_1 B$ if and only
if there are injective $m$-reductions from $A$
to $B$ and from $B$ to $A$.  
}

Both $\equiv_1$ and $\equiv_m$ are
$\Sigma^0_3\setminus \Pi^0_3$ equivalence relations on 
$\cantor$ (see \cite{rss}, Theorem 23).  
In Section \ref{1equivsec}, we will
show that both $\equiv_1$ and $\equiv_m $ 
are $\Pi^0_2$-graphable with diameter $2$.

(D) So far, we have only mentioned equivalence
relations on uncountable Polish spaces, but arithmetical
equivalence relations on $\N$ are also of interest for
these graphability questions.  For example, consider the equivalence
relation $E$ on $\N$ defined by
\[
n Em \iffdef W_n \equiv_1 W_m,
\]
where (as usual) $W_e$ is the domain of the $e$th
computable partial function on $\N$.  
This is a $\Sigma^0_3$-complete equivalence relation on $\N$
(see \cite{fokina2012equivalence}).
In Section \ref{indexsection}, it will be shown
that $E$ is $\Pi^0_2$-graphable with diameter $2$; in fact, this
follows from a more general result on arithmetical index relations,
see Theorem \ref{indexer}.

(E) Let $\mathcal{L}$ be a countable language,
and let $X_\mathcal{L}$ be the space of $\mathcal{L}$-structures on $\N$.
Isomorphism of $\mathcal{L}$-structures, $\cong_\mathcal{L}$,
is an 
analytic equivalence relation on $X_\mathcal{L}$.  
\cite{akl} proves that for any countable
language $\mathcal{L}$, the isomorphism equivalence relation
$\cong_\mathcal{L}$ is Borel graphable.

If we impose computability requirement on the isomorphisms,
then we get the arithmetical equivalence relation
\[
x \ciso y \iffdef \text{$x$ and $y$ are computably isomorphic}.
\]
When $\mathcal{L}$ is a computable language, 
$\ciso$ is a $\Sigma^0_3$ equivalence relation.
In Section \ref{computableisos}, 
we will show that for any computable language,
$\ciso$ is $\Pi^0_2$-graphable with diameter $2$.
The result and its method of proof will extend to many
related arithmetical equivalence relations, e.g., computable
isomorphisms of linear orders and computable
biembeddability of $\mathcal{L}$-structures.

(F) If $E$ is an equivalence relation on $X$, 
then the \textbf{Friedman-Stanley
jump} $E$, denoted $E^+$, is the equivalence relation
on $X^\N$ defined by
\[
(x_i)E^+(y_i) \iffdef \{[x_i]_E : i\in \N\} = \{[y_i]_E : i\in \N\}.
\]
It is easy to see that if $E$ is arithmetical, $E^+$ is
also arithmetical.  We will study
graphings for Friedman-Stanley jumps in Section \ref{fsjumps}
and show that from a graphing of $E$ of finite diameter, we can
arithmetically define a graphing of $E^+$.

\smallskip

\textbf{Notation and conventions.} 
Variables $i, j, k, n, m, \ell$, etc., will range over
natural numbers, and $x, y, z,$ etc., will range over
elements of
(uncountable) Polish spaces.  In definitions, quantifiers
like $(\forall n)$, $(\exists m)$, etc.,  will be understood to range
over natural numbers, where $(\forall x)$, $(\exists y)$, etc.,
range over elements of whatever Polish space is currently
under consideration.

For $e\in \N$ and $x\in \cantor$, $\varphi_e^x:\N\rightarrow \N$ denotes
the $e$th  Turing machine run on oracle $x$. 
$\varphi_e$ is the $e$th Turing machine
run on the oracle of all zeros.  

If $\sigma$ is a finite sequence and $x\in 2^\N$, then
$\sigma^\wedge x\in \cantor$ is the concatenation of $\sigma$ followed
by $x$. 

As usual, we will often conflate subsets of $\N$ with 
    elements of $\cantor$, and for $A\subset \N$ we will use
    $A(n)$ to denote the value of the characteristic function
    of $A$ at $n\in \N$.  

For a binary relation $R$ on $X$, we will use both
$R(x, y)$ and $xRy$ to denote that $x$ is $R$-related
to $y$.  Moreover, for $A\subset X\times Y$ and $x\in X$,
we will frequently use the notation $A_x:=\{y\in Y : A(x, y)\}$
for the $x$-section of $A$.  

We work with effective pointclasses, so our setting will be 
recursive Polish spaces and this is what we mean 
whenever we say that $X$ is a space.\footnote{
See \cite{rps} for details about the 
definition of recursive Polish spaces.
}
However, most of our equivalence relations are on 
$\N$, $2^\N$, the Baire space $\N^\N$, and
spaces which are computably isomorphic to 
products of these spaces, 
so one can safely restrict to such
spaces without losing too much generality.

We will use the notation $\forall^\N \Gamma$ to denote
the pointclass obtained by placing
a universal quantifier over $\N$ in front of $\Gamma$ 
relations.
For example, $\forall^\N \Sigma^0_2 = \Pi^0_3$ and
$\forall^\N \Pi^0_4 = \Pi^0_4$.  We define
pointclasses $\exists^\N \Gamma$, 
$\forall^\N\exists^\N \Gamma$,
etc., in a similar way.

Several of the equivalence relations we consider, notably
$E_0$ and $\equiv_T$, can be interpreted as relations on 
$\cantor$ or $\N^\N$.  By default, we consider them to be
equivalence relations on $\cantor$, unless we 
specify otherwise.

\smallskip

\textbf{Acknowledgments.}  I would like to thank
Anton Bernshteyn, Alexander Kechris, and Andrew Marks 
for many helpful
discussions on these topics.  
I am very grateful to Forte Shinko and Felix Weilacher
for proving and communicating Theorem \ref{shinko},
which establishes the optimality of many of the graphing
results to follow.  
Finally, thanks
to Patrick Lutz for many helpful comments on early
versions of proofs, in particular a comment that
strengthened the result of Theorem \ref{1equiv} on 
$1$-equivalence.

\section{Basic examples and properties}
\label{basic}

On any recursive Polish space $X$, equality
on $X$, denoted $=_X$, is a
$\Pi^0_1$ equivalence relation.  It is computably graphable,
using the trivial graph $G=\emptyset$.
More generally, let $\Gamma$ be a 
pointclass that contains $\Sigma^0_1$.
If $E$ is a $\Gamma$ equivalence relation on $E$, then
it has a $\Gamma$ graphing, $G:= E\setminus (=_X)$.  Of course, on its own this observation
is not interesting since 
we are interested in graphings of $E$ that have a simpler definition
in the arithmetical hierarchy than $E$ has.\footnote{
These first two observations are trivial, but they do have their uses;
see Section \ref{fsjumps} on Friedman-Stanley jumps.
}

Almost all of the proofs of graphability with diameter $2$ results 
will make use of the following simple lemma.

\begin{lemma}
\label{obvious}
    Let $\Gamma$ be a pointclass that contains $\Sigma^0_1$
    and is closed under $\&$ and $\vee$.  Let $E$ be an equivalence
    relation on a space $X$.  If there exists a binary relation
    $R$ on $X$ which is in $\Gamma$ and satisfies
    \begin{enumerate}[(i)]
    \item for all $x, y\in X$, $xRy$ implies $xEy$; and
    \item for all distinct $x, y\in X$ which are $E$-equivalent,
    there exists $z\in X$ such that $xRz$ and $yRz$,
    \end{enumerate}
    then $E$ is $\Gamma$-graphable with diameter $2$.  
\end{lemma}

\begin{proof}
    All we need to do is symmetrize $R$ and remove the diagonal.
    Formally, define a graph $G$ on $X$ by
    \[
    xGy \iffdef x\neq y \conj[xRy \vee yRx],
    \]
    and (easily) verify that $G$ is a $\Gamma$-graphing of $E$
    with diameter $2$.
\end{proof}

The proofs to follow will just define the relation $R$ and
implicitly use Lemma \ref{obvious}, leaving out the final
steps of symmetrizing and removing the diagonal.  In fact,
we will conflate $R$ with the final graphing 
by calling it $G$
from the start.

The next result establishes our first example
of a simpler graphing of an important arithmetical 
equivalence relation.  

\begin{proposition}[Folklore]
\label{turingequiv}
Turing equivalence $\equiv_T$ 
is $\Pi^0_2$-graphable with diameter 2.
\end{proposition}

\begin{proof}
    We first introduce some notation.  
    For any $n\in \N$, denote by $\sigma_n$ the
    length $n+1$ sequence of the form $(0, \dots, 0, 1)$,
    so that $\sigma_n(n)=1$ and $\sigma_n(i)=0$ for $i<n$.
    Also, fix $e^*$ so that $\varphi_{e^*}^x= x$ for all
    $x\in \cantor$.
    
    \ignore{For every $n\in \N$, there is a program $e_n$
    which satisfies $\varphi_{e_n}^x=\sigma_n{}^\wedge x$
    for every $x\in\cantor$, and we can find $e_n$ computably.
    Also, fix a program $e^*$ such that for every $n\in \N$,
    \[
    \varphi_{e^*}^{\sigma_n{}^\wedge x} = x.
    \]
    Let $g(n)$ be a 
    computable function such that $e^*, n, e_n<g(n)$
    for every $n$.}
    
    Define $G(x, z)$ to hold exactly when 
    $z$ is of the form $\sigma_n{}^\wedge y$ and there
    exists $e_0, e_1<n$ such that
    $\varphi_{e_0}^x=y$ and $\varphi_{e_1}^y=x$.   

    It is clear that $G$ is $\Pi^0_2$ and that
    $G(x, z)$ implies $x\equiv_T z$.  We claim
    that for any distinct $x\equiv_T y$,
    there is $z$ such that $G(x, z)$ and $G(y, z)$
    both hold.  
    For such $x, y$, 
    pick $n$
    large enough so that $e^*<n$ and
    there are $e_0, e_1<n$ such that
    $\varphi_{e_0}^x=y$ and $\varphi_{e_1}^y=x$.
    Then we take $z:= \sigma_n{}^\wedge y$.
    It is readily verified that $G(x, z)$
    and $G(y, z)$ both hold. Note that the fact that
    $G(y, z)$ holds uses that $e^*<n$.
\end{proof}

Examples of equivalence relations on $\N$
that admit simpler graphings can be obtained by
relativizing the following result from \cite{akl}.

\begin{theorem}[\cite{akl}] 
\label{aklceer}
    Suppose $E$ is a $\Sigma^0_1$ equivalence relation on $\N$,
    all of whose equivalence classes are infinite.  
    Then, $E$ is computably graphable with diameter $2$.
    \end{theorem}

    \begin{corollary}
    \label{ceercor}
    Suppose $E$ is a $\Sigma^0_{n}$ equivalence relation on $\N$,
    all of whose equivalence classes are infinite.  
    Then, $E$ is $\Delta^0_n$-graphable with diameter $2$.
    \end{corollary}

    \begin{proof}
        For $m\in \N$, let $\emptyset^{(m)}$ denote the $m$th Turing jump
        of $\emptyset$. 
        Since $E$ is $\Sigma^0_n$, it is $\Sigma^0_1(\emptyset^{(n-1)})$.
        By the relativized version of Theorem \ref{aklceer},
        $E$ has a $\emptyset^{(n-1)}$-computable graphing $G$ of diameter $2$.
        We are done since $\emptyset^{(n-1)}$-computable relations are exactly the $\Delta^0_n$ relations.
    \end{proof}

The rest of this section will be dedicated to some
closure properties of graphability, the first being about
computable
reducibility.  If $E$ is an equivalence relation
on $X$ and $F$ is an equivalence relation on $Y$, then a 
\textbf{computable reduction} of $E$ to $F$ is
a computable map $f:X\rightarrow Y$ such that 
$xEx'$ if and only if $f(x)Ff(x')$ for all $x, x'\in X$.
The computable reduction $f$ is said to be \textbf{invariant}
if its range is a union of $F$-classes, i.e.,
whenever $y\in Y$ is $F$-equivalent to some element of
the range of $f$, $y$ is also in the range of $f$.

\begin{proposition}
\label{invariantprop}
Let $\Gamma$ be a pointclass which contains both
$\Sigma^0_1$ and
$\Pi^0_1$, and is closed under $\vee$, $\&$, and
computable substitutions. 
Let $E$ and $F$ be equivalence relations on 
spaces $X$ and $Y$, respectively.   
If $E$ is invariantly 
computably reducible to $F$ and $F$ is $\Gamma$-graphable,
then $E$ is $\Gamma$-graphable.  Moreover, if $F$
has finite diameter $\ell$, then $E$ has finite diameter
$\leq \ell$.  
\end{proposition}

\begin{proof}
    Let $f:X\rightarrow Y$ be an invariant computable reduction
    of $E$ to $F$ and let $G$ be a $\Gamma$-graphing of $F$.
    Define a graph $H$ on $X$ by
    \[
    xHx' \iffdef x\neq x' \conj [f(x)=f(x') \vee f(x)Gf(x')]
    \]
    $H$ is in $\Gamma$ by our assumptions on $\Gamma$. 
    It is easy to check that $xHx'$ implies $xEx'$, so
    what is left to check is that any distinct $E$-equivalent
    element are connected by a path in $H$.  Suppose
    $x, x'\in X$ are distinct and $E$-equivalent.  
    If $f(x)=f(x')$, then we have $xHx'$, so we assume $f(x)\neq f(x')$.
    Since $f$ is a reduction, we have $f(x)Ff(x')$, hence 
    there exists a  $G$-path 
    $f(x), y_0, \dots, y_{k-1}, f(x')$. 
    Using that $f$ is invariant, we can find 
    $x_0, \dots, x_{k-1}\in X$ with $f(x_i) = y_i$ for all $i<k$.  Thus,
    $x, x_0, \dots, x_{k-1}, x'$ is a path in $H$.  
\end{proof}

For a pointclass $\Gamma$, an equivalence relation
$E$ on $\N$ is said to be \textbf{universal} under
computable reductions for $\Gamma$
equivalence relations on $\N$ if every $\Gamma$ equivalence
relation on $\N$ computably reduces to $E$.\footnote{
Since we are most concerned with lightface pointclasses,
we will always consider universality under \textit{computable} reductions.
}

Before we state our next result, we will quickly give the construction
of a universal $\Sigma^0_n$ equivalence relation on $\N$ for $n\geq 1$.
Let $U\subset \N\times \N^2$ be universal for $\Sigma^0_n$ 
subsets of $\N^2$, i.e., $U$ is $\Sigma^0_n$ and
its sections $U_i$ are exactly the $\Sigma^0_n$ subsets of $\N^2$.  
Define $V\subset \N\times \N^2$ so
that each $V_i$ is obtained from $U_i$ by symmetrizing, taking
the transitive closure, and adding the diagonal. 
Thus, each $V_i$ is an equivalence relation.
It is easy to check that $V$ is again $\Sigma^0_n$.
Moreover,
for every $\Sigma^0_n$ equivalence relation $E$ on $\N$, $E=U_i$ for some $i\in \N$. Since $U_i$ is already
an equivalence relation, $E=U_i=V_i$.  
Now, define the equivalence relation 
$F\subset (\N\times \N)^2$ by setting
\[
(i, n)F(j, m) \iffdef i=j \conj E(i, n, m).
\]
By the previously mentioned property of $E$,
every $\Sigma^0_n$ equivalence relation reduces to $F$.
You can use a computable isomorphism between $\N$ and $\N^2$
to make $F$ an equivalence relation on $\N$, which
gives a universal $\Sigma^0_n$ equivalence relation on $\N$.

\begin{proposition}
    For every $n\geq 1$, there is a universal $\Sigma^0_n$
    equivalence relation on $\N$ which is $\Pi^0_{n-1}$-graphable with diameter $3$.
\end{proposition}

\begin{proof}
    Let $E\subset \N^2$ be a universal $\Sigma^0_n$ equivalence
    relation on $\N$.  Pick a $\Pi^0_{n-1}$ relation $R\subset \N^3$
    such that 
    \[
    nEm \iff (\exists k)R(n,m, k).
    \]
    Now, define a graph $G$ on $\N\cup \N^3$ by only putting
    edges between $a\in \N$ and $(n, m, k)$ when $a\in \{n, m\}$ 
    and $R(n, m, k)$ holds.  $G$ is easily a $\Pi^0_{n-1}$ graph.
    Let $E_G$ be its connectedness relation.

    We show that $G$ has diameter $3$.  There are 
    two cases to consider.  First, suppose
    $nE_Gm$. Then, there is a $G$-path of the form
    \[
    n, (n, m_0, k_0), m_0, (m_0, m_1, k_1), 
    \dots, m_{\ell-1}, (m_{\ell-1}, m, k_\ell), m.
    \]
    Each tuple is in $R$, so $nEm_0Em_1\cdots m_{\ell-1}Em$.  Thus, $nEm$.
    But then there is $k$ with $R(n, m, k)$, which means $n, (n, m, k), m$
    is a path in $G$.  The other case is if $nE_G(m, \ell, k)$.  
    A similar argument shows that $nEm$.  Thus, there is a $k'$
    with $R(n, m, k')$, so that $n, (n, m, k'), m, (m, \ell, k)$
    is a path in $G$.  

    We now show that $E_G$ is universal for $\Sigma^0_n$
    equivalence relations on $\N$.  Let $F$ be a $\Sigma^0_n$
    equivalence relation on $\N$.  By universality of $E$,
    there is a computable reduction $f:\N\rightarrow \N$
    from $F$ to $E$.  Consider $f$ now as a function from $\N$
    to $\N\cup \N^3$.  We claim that $f$ is a reduction from
    $F$ to $E_G$.  Indeed, $nFm$ iff $f(n)Ef(m)$ iff there
    is $k$ such that $R(f(n), f(m), k)$ iff there is $k$ such that
    $f(n), (f(n), f(m), k),
    f(m)$ is a path in $G$ iff $f(n)E_Gf(m)$.  Note that
    the last ``iff'' uses that any $E_G$-equivalent numbers
    have a path of length 2 between them, which was established in
    the previous paragraph.
    \ignore{
    Suppose $nEm$.  Then, $f(n)Ef(m)$ so that 
    there is $k\in \N$ with $R(f(n), f(m), k)$.  Thus,
    $f(n), (f(n), f(m), k), f(m)$ is a path in $G$, hence
    $f(n)E_Gf(m)$.  Conversely, suppose
    $f(n)E_Gf(m)$. Then, there is a $G$ path of the form
    \[
    f(n), (f(n), m_0, k_0), m_0, (m_0, m_1, k_1), 
    \dots, (m_{\ell-1}, f(m), k_\ell), f(m).
    \]
    Each tuple is in $R$, so $f(n)Em_0Em_1\cdots m_{\ell-1}Ef(m)$.  Thus, $f(n)Ef(m)$, from
    which is follows $nFm$.}
\end{proof}

In the next section, Proposition \ref{sigmanotpi}
will establish (in part) that for $n\geq 1$, there
are $\Sigma^0_n$ equivalence relations on $\N$
which are not $\Pi^0_n$-graphable.  Since for $n\geq 1$
there are universal $\Sigma^0_n$ equivalence relations on $\N$ which
are $\Pi^0_{n-1}$-graphable, we have the following result.

\begin{corollary}
    For every $n\geq 1$, $\Pi^0_n$-graphability is not
    closed under computable (non-invariant) reduction.
\end{corollary}

We end this section by pointing out a few
closure properties about products of equivalence
relations.

\begin{proposition}
Let $\Gamma$ be a pointclass which contains $\Sigma^0_1$ and
$\Pi^0_1$, and is closed under
$\vee$, $\&$, and computable substitutions.
    \begin{enumerate}[(i)]
    \item Let $E$ and $F$ be equivalence relations
    on spaces $X$ and $Y$, respectively. Let $E\times F$
    be the product equivalence relation on $X\times Y$ defined
    by
    \[
    (x, y) E\times F (x', y') \iffdef xEx' \conj yFy'.
    \]
    If $E$ and $F$ are both $\Gamma$-graphable,
    then $E\times F$ is $\Gamma$-graphable.  Moreover,
    if $E$ and $F$ are $\Gamma$-graphable with finite
    diameters $k$ and $\ell$, respectively, then 
    $E\times F$ is $\Gamma$-graphable with diameter
    $\max(k, \ell)$.  

    \item For each $i\in \N$, let $E_i$ be an equivalence
    relation on space $X_i$, and suppose $X=\prod_i X_i$
    is also a recursive Polish space.  Suppose there exists
    finite diameter
    graphs $G_i$ on $X_i$ such that each $G_i$ is
    a graphing of $E_i$ and the $G_i$ are $\Gamma$ uniformly in $i$.
    If the diameters of $G_i$ are uniformly bounded, then 
    $\prod_i E_i$ is $\forall^\N\Gamma$-graphable with diameter equal
    to the maximum of the diameters of $G_i$.  In particular,
    if $E$ is $\Gamma$-graphable with diameter $k$, then 
    the infinite product $\prod_i E$ is $\forall^\N\Gamma$-graphable
    with diameter $k$.  
    \end{enumerate}
\end{proposition}

\begin{proof}
    For (i), let $G_E$ and $G_F$ be $\Gamma$-graphings of
    $E$ and $F$, respectively.  Then define a binary relation
    $H$
    on $X\times Y$ by
    \[
    (x, y) H (x', y') \iffdef [x=x' \vee xG_Ex'] \conj
    [y=y' \conj yG_F y'].
    \]
    The desired $\Gamma$ graphing can be obtained by
    removing the diagonal from $H$.  

    The proof of (ii) is quite similar.  We define
    \[
    (x_i)H(x_i') \iffdef (\forall i)[x_i=x_i'\vee x_iG_ix_i'].
    \]
    The assumption that the diameters of the $G_i$ are uniformly
    bounded is used to ensure that there is a finite path
    between any $\prod_i E_i$-equivalent sequences.  The shortest
    such path will of course have length no larger than the
    largest diameter of the $G_i$.  
\end{proof}

\ignore{
As mentioned before, the equivalence relation $E_0$ on $\cantor$
has an extremely important graphing, $G_0$.
$G_0$ is $\Sigma^0_2$, the same complexity as $E_0$.  

\begin{proposition}
    $E_0$ is $\Pi^0_1$-graphable (with diameter 2).
\end{proposition}
}

\section{Negative graphability results}
\label{negsection}

In this section, we will establish several results about
when graphings with certain types of
definitions are not possible. 

\begin{proposition}
    Let $n\geq 1$.
    If $E$ is an equivalence relation on $\N$ which 
    is not $\Sigma^0_n$, then $E$
    is not $\Sigma^0_n$-graphable.  In particular,
    for every $n\geq 1$ there are $\Pi^0_n$ equivalence
    relations on $\N$ which are not $\Sigma^0_n$-graphable.
\end{proposition}

\begin{proof}
    Suppose towards a contradiction that $G$
    is a $\Sigma^0_n$-graphing of $E$.  Then,
    \[
    nEm \iff (\exists k_0, \dots, k_{\ell-1})
    [nGk_0Gk_1\cdots k_{\ell-1}Gm].
    \]
    Using standard sequence coding techniques,
    this shows that $E$ is $\Sigma^0_n$, which is a 
    contradiction.
\end{proof}

On uncountable spaces, we also have $\Pi^0_n$ equivalence
relations with no simpler graphings.  We will see in the
sequel that many of our graphings take advantage of having
infinite equivalence classes, so the next result also
points out that having infinite equivalence classes
is in general not enough to guarantee simpler graphings.  

\begin{proposition}
    Let $n\geq 1$.
    There is a $\Pi^0_n$ equivalence relation on $\N^\N$,
    all of whose classes are countably infinite, which is
    not $\Sigma^0_n$-graphable.
\end{proposition}

\begin{proof}
    We work on the space $\N^2\times \N^\N$, which is 
    computably isomorphic to $\N^\N$.  Pick some $A\subset \N^\N$
    which is $\Pi^0_n\setminus \Sigma^0_n$.  Now, define
    an equivalence relation $E$ by 
    \[
    (n, m, x)E(n', m', x') \iffdef [m=m' \conj x=x'] \vee [ x=x'\in A].
    \]
    It is clear that $E$ is $\Pi^0_n$ and that all equivalence classes
    are countably infinite. 
    
    Suppose towards a contradiction that $G$ is a $\Sigma^0_n$
    graphing of $E$.  We claim that
    \[
    x\in A \iff (\exists n, n', m, m')[m\neq m' \conj 
    (n, m, x)G(n', m', x)],
    \]
    which is enough since the claim implies $A$ is $\Sigma^0_n$.
    The right-to-left direction is immediate.  Suppose now
    that $x\in A$ but the right-hand-side fails.  Then, every
    $G$-adjacent point to 
    $(n, 0, x)$ is of the
    form $(n', 0, x)$ for some $n'\neq n$.  A simple induction
    shows that the $G$-connected component of $(0, 0, x)$ does not
    contain, say, $(0, 1, x)$, which is $E$-equivalent to $(0, 0, x)$.
    Thus, $G$ is not a graphing of $E$.  
\end{proof}

We can use a similar strategy to produce a $\Sigma^0_n$
equivalence relation which does not have any simpler
graphings.  In this case, the equivalence relation
will have finite equivalence classes.

\begin{proposition}
\label{sigmanotpi}
For any $n\geq 2$, there exists 
a $\Sigma^0_n$ equivalence
relation on $2^\N$ which is not $\Pi^0_n$-graphable.  
The same result also holds for the case where $n\geq 1$ and
the space is $\N$.  
\end{proposition}

\begin{proof}
Fix $n\geq 2$ and fix some $A\subset 2^\N$ which is $\Sigma^0_n\setminus \Pi^0_n$.
Define $E$ on $\{0, 1\}\times 2^\N$ by
\[
(i, x) E(j, y) \iffdef (i=j \conj x=y) \vee (x=y\in A).
\]
Easily, $E$ is a $\Sigma^0_n$ equivalence relation.  
Suppose towards a contradiction that 
$G$ is a $\Pi^0_n$ graphing of $E$.  Then,
\[
x\in A \iff (0, x) G (1, x),
\]
contradicting that $A$ is not $\Pi^0_n$.  

The claim about $\N$ is proved identically, with the
case $n=1$ also working because
equality on $\N$ is computable.
\end{proof}

 We know that there are $\Pi^0_n$ equivalence relations
    with infinite classes that have no simpler graphings and
    that, for $\Sigma^0_n$ equivalence relations on $\N$, 
    infinite classes are enough to guarantee a 
    $\Delta^0_n$ graphing. This leads us to the following
    open problem.
    
\begin{question}
    Is there a $\Sigma^0_n$ equivalence relation 
    on an uncountable space, all of whose classes are infinite,
    which is not $\Pi^0_n$-graphable? 
\end{question}

As previously mentioned, many of the graphability results
to follow produce graphings of diameter $2$.  Next,
we will see that there are indeed situations in which
finite diameter graphings are not possible.

\begin{theorem} 
\label{nonfinitediam}
Let $X$ be a compact space and
let $E$ be an equivalence relation on $X$.  If
$E$ has an equivalence class which is not closed,
then $E$ does not have any finite diameter closed graphings.
\end{theorem}

\begin{proof}
    Suppose towards a contradiction that $G$ is a
    closed graphing of $E$ with finite diameter $k$.  
    Let $C$ be an $E$-class which is not closed, and fix $x\in C$.
    Define a $k$-ary relation $R$ by
    \[
    R(x_1, \dots, x_k) \iffdef xGx_1 \conj (\forall i<k)[x_i=x_{i+1}
    \vee x_i Gx_{i+1}].
    \]
    Clearly, $R$ is closed.  Since $G$ is a graphing of $E$,
    $R(x_1, \dots, x_k)$ implies that
    $x_1, \dots, x_k\in C$.
    Moreover, since the diameter of $G$ is $k$, 
    for every $y\in C$ there exist 
    $x_1, \dots x_{k-1}\in C$ such that
    $R(x_1, \dots, x_{k-1}, y)$. 

    Since $C$ is not closed, we may fix some $y\notin C$ which 
    is in the closure of $C$.  Then,  
    we can find some sequence  $(y_n)$ in $C$ which converges
    to $y$.  For every $n$, pick $x_{1, n} \dots, x_{k-1, n}\in C$
    such that $R(x_{1, n}, \dots, x_{k-1, n}, y_n)$ holds.
    Using the compactness of $X$, we may pass to a converging
    subsequence, so we may assume that 
    $(x_{1, n}, \dots, x_{k-1, n}, y_n)$
    converges to some $(x_1, \dots, x_{k-1}, y)$.  Since $R$
    is closed, it follows that $R(x_1, \dots, x_{k-1}, y)$ holds,
    which implies the contradiction that $y\in C$.
\end{proof}

\begin{corollary}
    There is a $\Sigma^0_2$ equivalence relation
    which is $\Pi^0_1$-graphable with infinite diameter,
    but does not have any closed graphings of finite diameter.
\end{corollary}

\begin{proof}
Consider the orbit equivalence relation $E$ of
an irrational rotation of the circle.  If the irrational
angle $\theta$ is a computable real, then $E$ is $\Sigma^0_2$.  
The graph $G$ with edges between $x, y$ when they are one
$\theta$-rotation apart is a $\Pi^0_1$-graphing with 
infinite diameter.  It is well known that every class of $E$
is dense, hence not closed.  
By Theorem \ref{nonfinitediam},
$E$ has no closed graphings of finite diameter.
\end{proof}

By a result of Clemens \cite{clemens}, $E_0$ 
can be generated by a single homeomorphism
on $2^\N$, which in particular gives us a closed graphing
of $E_0$ with infinite diameter.  Of course,
Theorem \ref{nonfinitediam} also applies to $E_0$.
Thus, $E_0$ is another example
of a $\Sigma^0_2$ equivalence relation that has a closed
graphing of infinite diameter but no closed graphings of finite
diameter.  For finite diameter graphings of $E_0$, 
the following result is the
best we can do.

\begin{proposition}
    $E_0$ is $\Delta^0_2$-graphable with diameter $2$.
\end{proposition}

\begin{proof}
    Define a binary relation $G$ 
    so that $G(x, y)$ holds exactly when 
    $y$ begins with some $\sigma_n$ (as
    defined in the proof of Proposition \ref{turingequiv}),
    and $x(i)=y(i)$ for all $i\geq n+1$.  Clearly,
    $G$ is $\Pi^0_2$ and $xGy$ implies $xE_0 y$. 
    Given $xE_0 y$,
    pick some $n$ so that $x(i)=y(i)$ for all $i>n$,
    then define $z$ so that $\sigma_n\subset z$
    and $z(i) = x(i)$ for all $i\geq n+1$.  Then,
    it is easy to check that $xGz$ and $yGz$.  
\end{proof}

The compactness of the underlying space in Theorem
\ref{nonfinitediam} is essential.  To see this,
we will consider the equivalence relation $E_0(\N)$ of
eventual equality on the Baire space $\N^\N$ (which, of
course, is not compact).  Note
that all of the classes of $E_0(\N)$ are dense, so in particular
not closed.

\begin{proposition}
   $E_0(\N)$ is $\Pi^0_1$-graphable with diameter $2$.
\end{proposition}

\begin{proof}
    Define a graph $G$ on $\N^\N$ to have an edge
    between $x$ and $y$ when $x(0)\neq y(0)$ and $x(i)=y(i)$
    for all $i> \max(x(0), y(0)))$.   
    Clearly, $G$ is $\Pi^0_1$ and $xGy$ implies 
    $x$ and $y$ are eventually equal.  
    Let $x$ and $y$ be $E_0(\N)$-equivalent.  We must
    find $z$ so that the pairs $x, z$ and $z, y$ are $G$-adjacent.
    To do this, pick $k>x(0), y(0)$ so that $x(i)=y(i)$ for all $i>k$.  
    Then, define $z$ so that $z(0) = k$ and
    $z(i) = x(i)$ for $i>0$.  It is routine to check
    that $xG z$ and $zGy$ both hold.  
\end{proof}

Another example of a $\Sigma^0_2$ equivalence
relation on a non-compact space which has a closed graphing
of diameter $2$ is the Vitali equivalence relation, denoted $E_V$.  
$E_V$ is the equivalence relation on $\R$ given by
\[
xE_Vy \iffdef x-y \in \mathbb{Q}.
\]
Note that the equivalence classes of $E_V$ are all dense.

\begin{proposition}
    The Vitali equivalence relation $E_V$ is closed graphable
    with diameter $2$.
\end{proposition}

\begin{proof}
    Fix an enumeration $q_1, q_2, \dots$ of 
    $\Q\setminus \Z$. 
    Define a binary relation $G$ on $\R$ 
    so that $G(x, y)$ holds when
    \begin{enumerate}[(i)]
    \item $x-y\in \Z\setminus\{0\}$; or
    \item $x\leq y$, $1\leq y$, and there is a positive integer
    $n\leq y$ such that $y-x-q_n\in \Z$.  
    \end{enumerate}
    Immediately, $xGy$ implies $xE_Vy$.  
    Also, $G$ is irreflexive.\footnote{
    Note that this is important since we want a closed graphing.  $=_\R$ is closed, so removing it from
    a closed relation would make the resulting relation
    $\Delta^0_2$.}
    To see this, note
    that $x-x-q_n\notin \Z$ since $q_n$ is not
    an integer.  
    
    We will next show that the symmetrization of $G$ 
    is a graphing of $E_V$.  
    Suppose $x< y$ are $E_V$-equivalent. If $x-y\in \Z$, then
    $G(x, y)$ holds by virtue of $(i)$, so assume $x-y\notin \Z$.
    Fix $n$ with $y-x = q_n$.  Pick $k\in \N\setminus\{0\}$
    large enough so that
    $y+k$ is larger than $n$.  Clearly,
    $yG(y+k)$ by clause (i).  Moreover, $xG(y+k)$ because
    $(y+k)-x = (y-x)+k = q_n +k$, where $n\leq y+k$.  

    What is left to show is that $G$ is closed.  Clause (i)
    is clearly a closed condition, 
    so we only have to show that (ii) is a closed condition.  
    Suppose $(x_n)$ and $(y_n)$ are sequences in $\R$
    converging to $x$ and $y$, respectively, and that for
    each $n\in \N$, the pair $(x_n, y_n)$ satisfies (ii).
    Clearly, $x\leq y$ and $1\leq y$.  Since $y_n\rightarrow y$,
    for sufficiently large $n$ we have
    $y_n-x_n-q_{i_n}\in \Z$ for some $1\leq i_n\leq \lfloor y \rfloor$.  
    By passing to a subsequence, we may assume there is
    a fixed positive integer $i_0\leq \lfloor y \rfloor$ such that $y_n-x_n-q_{i_0}\in \Z$
    for all $n$.  Since $Z$ is closed, $y-x-q_{i_0}\in \Z$ and $xGy$ holds.  
\end{proof}

The following non-graphability result will establish the 
optimality of many of the positive results to follow.
It was proved and communicated by Forte Shinko and 
Felix Weilacher.

\begin{theorem}[Shinko-Weilacher]
\label{shinko}
Let $E$ be an equivalence relation on $\cantor$.
$E$ is $\Sigma^0_2$-graphable if and only if $E$
is $\Sigma^0_2$.  Moreover, this statement relativizes to any
real parameter.
\end{theorem}

\begin{proof}
    The right-to-left direction is 
    immediate using $G:=E\setminus (=_{2^\N})$. 
    For the other direction, 
    suppose that $G$ is a $\Sigma^0_2$ graphing of $E$.
    Pick a $\Pi^0_1$ relation 
    $H\subset \N\times(\cantor)^2$ such that
    \[
    G(x, y) \iff (\exists k) H(k, x, y),
    \]
    so that $G$ is a union of the $\Pi^0_1$ graphs
    $H_k:=\{(x, y)\in (\cantor)^2 : H(k, x, y)\}$.
    
    Fix some computable bijection $F:2^\N\rightarrow (\cantor)^\N$. 
    We use the notation $(z)_i := F(z)(i)\in 2^\N$.
    Now, define
    \begin{multline*}
    R(n, m, x, y, z) \iffdef (z)_0=x \conj (z)_n=y \\
    \conj (\forall i<n)(\exists k<m)H(k, (z)_i, (z)_{i+1}).
    \end{multline*}
    It is clear that $R$ is $\Pi^0_1$ and that $R(n, m, x, y, z)$
    holds whenever $z$ codes a length $n$ 
    path from $x$ to $y$ which only
    uses edges from $H_k$ for $k<m$.  Using the compactness
    of $\cantor$, the condition $(\exists z)R(n, m, x, y, z)$
    is still $\Pi^0_1$.\footnote{
    Using K\"onig's lemma, the existence of such a $z$
    is equivalent to saying that some (computable) binary tree is infinite.}
    Thus, the equivalence
    \[
    xEy \iff (\exists n)(\exists m)(\exists z) R(n, m, x, y, z)
    \]
    shows that $E$ is $\Sigma^0_2$. 
\end{proof}

\begin{corollary}
\label{optimalitycor}
    None of the following equivalence relations are 
    $\Sigma^0_2$-graphable: $\equiv_T$, $\equiv_1$, $\equiv_m$,
    and $\cong_\mathcal{L}^c$, where $\mathcal{L}$ is a nontrivial
    computable relational language.
\end{corollary}

\begin{proof}
    These are all equivalence relations on $2^\N$ or, in the
    case of $\cong_\mathcal{L}^c$, on a space computably
    isomorphic to $2^\N$.  By Theorem \ref{shinko}, it is enough
    to show that these equivalence relations are not $\Sigma^0_2$.
    
    As mentioned previously, by results
    from \cite{rss}, 
    none of $\equiv_T$, $\equiv_1$, $\equiv_m$ are $\Pi^0_3$, so in 
    particular not $\Sigma^0_2$. It is also the
    case that $\cong_\mathcal{L}^c$ is not $\Pi^0_3$ when
    $\mathcal{L}$ is nontrivial, and it can be established by showing that $\equiv_1$ can be computably
    reduced to $\equiv_\mathcal{L}^c$.  Fix some relation symbol $R$
    in $\mathcal{L}$.  For a set $A\in 2^\N$, map it to the 
    $\mathcal{L}$-structure $x$ where the interpretation of
    every symbol other than $R$ is
    trivial, and where 
    \[
    R^x(a_0, \dots, a_{n-1}) \iff a_0=\cdots = a_{n-1}\in A.
    \]
    It is easily checked that this defines a computable reduction,
    so that $\cong_\mathcal{L}^c$ is also not $\Pi^0_3$.      
\end{proof}

In the sequel, we will establish that all of the equivalence
relation mentioned in Corollary \ref{optimalitycor} are
$\Pi^0_2$ graphable with diameter $2$, and so those are the simplest
possible in terms of arithmetical definability. 

Obviously, the proof of Theorem \ref{shinko} made crucial use
of the compactness of $\cantor$.  The situation for
$\Sigma^0_2$-graphings for equivalence relations on non-compact spaces
is left open.  The following is an interesting test case.

\begin{question}
    Let $\mathord{\equiv_T}(\N)$ be Turing equivalence on $\N^\N$.
    Is $\mathord{\equiv_T}(\N)$ $\Sigma^0_2$-graphable?  More generally,
    is there an equivalence relation on $\N^\N$ which is
    not $\Sigma^0_2$ but is $\Sigma^0_2$-graphable?
\end{question}


\section{1-equivalence of sets}
\label{1equivsec}

Recall that sets $A, B\subset \N$ 
are \textbf{$1$-equivalent}, written
$A\equiv_1 B$, if there is a computable bijection
$f:\N\rightarrow \N$ such that
$n\in A$ if and only if $f(n)\in B$
for all $n\in \N$.  Such an $f$ is called a \textbf{$1$-equivalence}.  Note that, given $A, B\subset \N$ and $e\in \N$,
checking
if $\varphi_e$ is a $1$-equivalence from $A$ to $B$
is $\Pi^0_2$, with the most complicated part being
checking that $\varphi_e$ is total.

In this section and in the sequel, we will make
use of the following notion.  
Call a bijection $f:\N\rightarrow \N$ a \textbf{finite swapping}
if it is the composition of finitely many transpositions.
Clearly, every finite swapping is a computable bijection.

    The following result will be generalized in Section
    \ref{computableisos} on computable isomorphisms. 
    Indeed, $1$-equivalence is just computable isomorphism for
    structures in the language $\mathcal{L}$ 
    with exactly one unary predicate.
    We present the proof for $1$-equivalence separately
    because this case will make the idea of 
    the proof of the more
    general theorem clearer and will also serve as special
    case in that argument.

    \begin{theorem}
    \label{1equiv}
    1-equivalence $\equiv_1$
    is $\Pi^0_2$-graphable with diameter 2.
    \end{theorem}

\begin{proof}
    We begin by defining a computable total
    $f:\N\rightarrow \N$
    which will provide computable bounds on searches for 
    programs of $1$-equivalences.  For $n\in \N$, let $f(n)$ be
    the maximum of all the program codes 
    of finite swapping functions which only use
    transpositions of numbers $<2n$ and of all the
    programs obtained from effectively composing such a bijection once with a program with some code $e<n$. 

    Now define a relation $G\subset \cantor\times \cantor$
    so that $G(A, B)$ holds exactly when
    \begin{enumerate}[(i)]
    \item $B$ is not $\emptyset$ or $\N$; and
    \item if $n$ is the
    least number $>0$ with $B(n-1)\neq B(n)$, then
    there is $e<f(n)$ such that $\varphi_e$ is
    an $1$-equivalence from $A$ to $B$. 
    \end{enumerate}
    It is a routine
    computation to check that $G$ is $\Pi^0_2$, and it is
    immediate that $G(A, B)$ implies $A\equiv_1 B$.  What
    is left to show is that for any pair of distinct $A, A'\subset \N$
    with $A\equiv_1 A'$, there exists $B\subset \N$ such
    that $G(A, B)$ and $G(A', B)$.  

    Let $A, A'\subset \N$ be distinct with $A\equiv_1 A'$.  
    Since $A, A'$ are distinct and $1$-equivalent,
    they are not equal to $\emptyset$ or $\N$.  Fix $e\in \N$
    such that $\varphi_e$ is an $1$-equivalence from $A'$ to $A$.

    Now, we will build the desired $B\subset \N$ by applying
    finitely many transpositions to $A$.  Fix $n\in \N$ such
    that $n>e$ and there are $m, m'<n$ with $m\in A$ and
    $m'\notin A$.  Among the numbers $<2n$, either at least $n$
    are in $A$ or at least $n$ are not in $A$.  Either way,
    we can swap them using finitely many transpositions so that
    the new set $B$ we obtain has $n$ as the least number
    with $B(n-1)\neq B(n)$.  Now, $A$ and $B$ are $1$-equivalent
    by a function generated by finitely many transpositions
    of numbers $<2n$, hence there is an $1$-equivalence between
    $A$ and $B$ whose index is $<f(n)$.  Thus, $G(A, B)$ holds.
    Moreover, $A'$ and $B$ are $1$-equivalent by a function
    obtained by composing $\varphi_e$, the $1$-equivalence from
    $A'$ to $A$, with the previously mentioned $1$-equivalence
    from $A$ to $B$.  
    Using that $e<n$ and our definition of $f(n)$,
    it follows that $G(A', B)$ also holds.  
\end{proof}

We can easily modify the proof of Theorem \ref{1equiv} to
get the following result about $m$-equivalence.

\begin{theorem}
    $m$-equivalence $\equiv_m$
    is $\Pi^0_2$-graphable with diameter 2.
\end{theorem}

\begin{proof}
    The proof is nearly identical.  The definition of
    $G$ should be altered slightly: 
    $G(A, B)$ holds exactly when
    $B$ is not $\emptyset$ or $\N$, and, if $n$ is the
    least number $>0$ with $B(n-1)\neq B(n)$, then
    there is $e,e'<f(n)$ such that $\varphi_e$ (respectively,
    $\varphi_{e'}$) is
    an $m$-reduction from $A$ to $B$ (resp., from $B$ to $A$).
    The proof proceeds as before, with the only substantive change
    being that in the construction of a set $B$ with
    $G(A, B)$ and $G(A', B)$, we choose $n$
    to be larger than both $e$ and $e'$, where these are fixed
    indices of $m$-reductions between $A$ and $A'$.
\end{proof}

\section{Index equivalence relations}
\label{indexsection}

Having examined the graphability of $\equiv_1$ and $\equiv_m$ as
relations on $\cantor$, we now address the graphability
of $1$-equivalence and $m$-equivalence 
as relations on c.e. indices. We will,
in fact, establish a result about any equivalence relation
on c.e. indices.  

An equivalence relation $E\subseteq \N\times \N$ is 
called an \textbf{index equivalence relation} if
whenever $\varphi_a=\varphi_{a'}$ and $\varphi_b=\varphi_{b'}$, we have
\[
aEb \iff a'Eb'.
\]
Recall that the Padding Lemma states that every computable partial
function $f$ is equal to $\varphi_e$ for infinitely many $e$.
It follows that every equivalence class
of an index equivalence relation is infinite.

Recall that Corollary \ref{ceercor} provides a 
nice upper bound on the
    definability of graphings for $\Sigma^0_n$ equivalence relations
    on $\N$, all of whose equivalence classes are infinite.
    The following theorem that shows we can do better
    than $\Delta^0_n$ graphings for index equivalence relations.

\begin{theorem}
\label{indexer}
    If $E$ is a $\Sigma^0_{n+1}$ index equivalence relation,
    then it is $\Pi^0_n$-graphable with diameter $2$.
\end{theorem}

\begin{proof}
    The Padding Lemma is
    true of any acceptable effective enumeration of
    the computable partial functions.  For concreteness,
    we will have
    in mind that indices code register machine
    programs.

    If $a, b, c\in \N$, we will say
    a program $e$ is of \textit{special form} (with respect
    to $a, b, c$) if $e$
    is the program
    \begin{center}
        ``add zero $a$ times to the 1st register, then add zero $b$ times to the 2nd register, then run program $c$.''
    \end{center}
    (Here, $a$ and $b$ are just thought of as numbers.)
    The relation
    \[
    \text{Form}(e) \iffdef e \ \text{is of special form with respect to some $a, b, c$}
    \]
    is clearly computable.  We mention some other obvious properties:
    \begin{enumerate}[(i)]
    \item For any $a, b, c\in \N$,
    there is a unique $e$ such that $e$ is of special form with
    respect to $a, b, c$.
    \item There are computable total functions 
    $\text{Code}_1, \text{Code}_2$, $\text{Main}$ such that whenever 
    $e$ is of special form with respect to $a, b, c$, then 
    \[
    \text{Code}_1(e)= a, \ \ \text{Code}_2(e)=b, \ \ \text{Main}(e)=c.
    \]
    \item If $e$ is of special form, then $\varphi_e = \varphi_{\text{Main}(e)}$.  Hence,  $\text{Main}(e)Ee$
    since $E$ is an index equivalence relation.
    \end{enumerate}
    
    Since $E$ is $\Sigma^0_{n+1}$, we can pick a $\Pi^0_n$ relation
    $R\subset \N^3$ such that
    \[
    aEb \iff (\exists n)~R(a, b, n).
    \]
    Now, define
    \begin{multline*}
    G(a, e) \iffdef \text{Form}(e) \\
    \conj \big [\text{Main}(e)=a
    ~\vee~ [\text{Code}_1(e)=a \conj R(\text{Main}(e), a, \text{Code}_2(e))],
    \end{multline*}
    which is clearly $\Pi^0_n$.  
    Next, we will show that if $G(a, e)$, then $aEe$.
    Suppose $G(a, e)$.  
    Then, $\text{Form}(e)$ and either $\text{Main}(e)=a$
    or we have $\text{Code}_1(e)=a$ and $R(\text{Main}(e), \text{Code}_1(e), \text{Code}_2(e))$ holds.
    In the former case, $a=\text{Main}(e)Ee$ by (iii) above.
    In the latter case,
     $R(\text{Main}(e), a, \text{Code}_2(e))$
    implies that $\text{Main}(e)Ea$. Thus,
    we have $eE\text{Main}(e)Ea$, hence $eEa$.  

    Now, assume $aEb$ and we will show that there exists 
    $e$ with $G(a, e)$ and $G(b, e)$. Pick $n$
    with $R(a, b, n)$.
    Now, let $e$ be of special form with 
    $\text{Main}(e) = a$, $\text{Code}_1(e)=b$ and $\text{Code}_2(e)=n$.
    It is immediately verified that $G(a, e)$ and $G(b, e)$
    both hold.  
\end{proof}

\begin{corollary}
    $1$-equivalence  and $m$-equivalence 
    of c.e. sets (as equivalence relations on indices) 
    are both $\Pi^0_2$-graphable.
\end{corollary}

\section{Computable isomorphisms}
\label{computableisos}

To begin, we will focus on the case of relational languages.
This case contains all of the interesting mathematics,
and the case of an arbitrary language will follow from it.
A countable relational language $\mathcal{L}$ is 
\textbf{computable} if there is an enumeration 
$\mathcal{L}=\{R_0, R_1, \dots\}$ such that the function
$i\mapsto \text{arity}(R_i)$ is computable.
Unless we specify otherwise, $\mathcal{L}$
will denote a computable
relational language, and we will denote its
relation symbols by $R_i$, $i\in \N$. 

The space of $\mathcal{L}$-structures on $\N$, denoted
$X_\mathcal{L}$, is a
recursive Polish space, usually realized as
\[
X_\mathcal{L} = \prod_{i}2^{(\N^{\text{arity}(R_i)})},
\]
which is computably isomorphic to $\cantor$. 
Note that  
\begin{equation}
\label{computstruc}
\{(x, i, \vec n)\in 
X_\mathcal{L}\times \N\times \N^{\text{arity}(R_i)}
: x\models R_i(\vec n)\}
\end{equation}
is computable.\footnote{
Note that, formally, we should use a sequence coding
so that the number of arguments in (\ref{computstruc}) does not change with 
the arities of $R_i$.  We will suppress the use of sequence
coding whenever it will not cause confusion.
}  
For $x\in X_\mathcal{L}$ and $R_i\in \mathcal{L}$,
we denote by $R_i^x$ the interpretation of $R_i$ 
in the structure $x$.

Let $x, y\in X_\mathcal{L}$.  An \textbf{$\mathcal{L}$-homomorphism}
from $x$ to $y$
is a total function $f:\N\rightarrow\N$ such that for every $R_i\in \mathcal{L}$,
\[
R_i^x(\vec a) \iff R_i^y(f(a_0), \dots, f(a_{k-1})) \qquad [k=\text{arity}(R)]
\]
If $f$ is an injection, then it is an \textbf{$\mathcal{L}$-embedding}.
If $f$ is a bijection, then $f$ is an \textbf{$\mathcal{L}$-isomorphism}.

Let $x\in X_\mathcal{L}$ and let $f:\N\rightarrow \N$ be a bijection.
The \textbf{pushforward of $x$ by $f$} is the structure 
$f_*x\in X_\mathcal{L}$ defined by
\[
R_i^{f_*x}(\vec a) \iffdef R_i^x(f^{-1}(a_0), \dots, f^{-1}(a_{k-1}))
\qquad [R_i\in\mathcal{L}, \ k=\text{arity}(R_i)]
\]
Clearly, $f$ is an isomorphism from $x$ to $f_*x$.  Moreover,
for any $x, y\in \mathcal{L}$, $f$ is an isomorphism from
$x$ to $y$ if and only if $y=f_*x$.

By a routine computation, 
``$\varphi_e$ is an isomorphism
from $x$ to $y$'' (as a relation of $e, x, y$) is $\Pi^0_2$.
Thus,
\[
x\ciso y \iff (\exists e)[\varphi_e \ \text{is an isomorphism
from $x$ to $y$}]
\]
is indeed $\Sigma^0_3$.

In an upcoming proof, we will construct a new structure $y$ from
some $x\in X_\mathcal{L}$ by describing finitely many 
transpositions
that create a finite swapping function $f$ and taking $y:=f_*x$.
It is often helpful to think of
numbers as labels on the elements
of $x$, rather than the elements themselves; 
this way, we are building
$y$ just by swapping around labels on the structure $x$.

For the case of $1$-equivalence of sets, 
we were able to store information in a set $B$
as the least number at which the characteristic function
of $B$
changes value.   Consider now the case of
a binary relation $R(n, m)$.  Supposing that the relation $R$
is non-trivial when 
restricted to pairs $(n, m)$ with $n\neq m$,
we could use a similar strategy, e.g., we
could store information as the least $n$
such that $R(0, n)$ fails.  However, it could be the case that, say,
$R(n, m)$ holds for all distinct $n, m$.  This is analgous in the case
of $1$-equivalence to the situation where $B=\N$; note that
this does not pose a problem for $1$-equivalence 
because $\N$ is only $1$-equivalent
to itself.  
However, in the case of the binary relation, 
there still may be nontrivial structure
in the diagonal relation $R(n, n)$.  But this is a unary relation,
so we could then apply our $1$-equivalence strategy 
to the diagonal relation.  What if the diagonal relation is
also trivial?  If so, then the whole binary relation $R(n, m)$
is extremely simple and easily dealt with.

The above discussion points us towards the idea that in the
case of a computable relational language $\mathcal{L}$ 
and $x\in X_\mathcal{L}$, we cannot
just look to the relations $R_i^x$; we need to also look
at the relations we can define from the $R_i^x$
using projection functions.  

For any $k>0$ and $i<k$, let $\pi^k_i:\N^k\rightarrow \N$
be the projection
\[
\pi^k_i(a_0, \dots, a_{k-1}) := a_i.
\]
If $\vec \pi=(\pi^k_{i(0)}, \dots, \pi^k_{i(n-1)})$ 
is a sequence of $k$-ary projections and $\vec a\in \N^k$,
we use the notation
\[
\vec \pi(\vec a) 
:= (\pi^k_{i(0)}(\vec a), \dots, \pi^k_{i(n-1)}(\vec a))
= (a_{i(0)}, \dots, a_{i(n-1)})
\]
We call such a $\vec \pi$ a \textbf{$k$-ary shuffling sequence
of length $n$}.

For an $n$-ary relation $R$ on $\N$ and $k\leq n$,
a \textbf{$k$-shuffle of $R$} is a $k$-ary relation $R_{\vec \pi}$
of the form
\[
R_{\vec \pi}(\vec a) \iffdef R(\vec \pi(\vec a)),
\]
where $\vec \pi$ is a $k$-ary shuffling sequence of length $n$.
For example, if $R$ is $3$-ary, then
\[
Q(a, b) \iffdef R(a, b, a)
\]
is a $2$-shuffle of $R$.

A shuffle of a shuffle of $R$ is again a shuffle of $R$.
More precisely, if $R$ is an $n$-ary relation, $\vec \pi$
is a length $n$ sequence of $k$-ary projections, and $\vec \rho$
is a length $k$ sequence of $\ell$-ary projections, then
$(R_{\vec\pi})_{\vec \rho}$
is an $\ell$-shuffle of $R$.  This follows from
the fact that if $\vec \pi = (\pi^k_{i(0)}, \dots, \pi^k_{i(n-1)})$
and $\vec \rho = (\pi^\ell_{j(0)}, \dots, \pi^\ell_{j(k-1)})$,
then 
\[
\vec \pi(\vec \rho(a_0, \dots, a_{\ell-1}))
= \vec \pi(a_{j(0)}, \dots, a_{j(k-1)}) 
=(a_{j(i(0))}, \dots, a_{j(i(n-1))}).
\]

A unary relation $R$ is a \textbf{coding relation}
if $\emptyset \subsetneq R \subsetneq \N$.
For $n>1$, an $n$-ary relation $R$ is a 
\textbf{coding relation}
if there is an injective sequence $\vec a$ of length $n-1$
such that
there exists $b, c$, distinct from all the $\vec a$, such that
$R(\vec a, b)$ and $\neg R(\vec a, c)$ hold.
$R$ has the \textbf{bad coding property} if
none of its shuffles are coding relations.

The next lemma shows us that if $R$ has the bad
coding property than it is very simple to define.
The \textbf{trivial language} is the language
with no non-logical symbols.  The only prime formulas
are $v=u$ (where $=$ is always interpreted as equality).

\begin{lemma}
    If $R$ has the bad coding property, then
    $R$ is definable with a (quantifier-free) formula
    in the trivial language.
\end{lemma}

\begin{proof}
    The proof is by induction on the arity $n$ of $R$.
    The case $n=1$ is trivial, since the assumption implies 
    that either
    $R=\N$ or $R=\emptyset$.

    Suppose now $n>1$ and the
    lemma holds for all arities less than $n$.
    First, we show that either $R(\vec a)$ holds for
    all injective $\vec a$, or $\neg R(\vec a)$ holds
    for all injective $\vec a$.  
    Suppose $R(\vec a)$ holds for some injective $\vec a$,
    and we will show that $R$ holds on all injective sequences.
    It is enough to show that for $i<n$ and 
    all $b\in \N$, different from
    the elements of $\vec a$, we have 
    \[
    R(a_0, \dots, a_{i-1}, b, a_{i+1}, \dots, a_{n-1}).
    \]
    Fix $i<n$.  Let 
    $\vec a_{\hat i}=(a_0, \dots, a_{i-1}, a_{i+1}, \dots, a_{n-1})$.
    Let $\vec \pi$ be the shuffling sequence with
    $\vec \pi (\vec a_{\hat i}, a_i)
    = \vec a$,
    so that $R_{\vec \pi}(\vec a_{\hat i}, a_i)$ holds.  Since $R_{\vec \pi}$ is not a coding relation, we must have that 
    $R_{\vec \pi}(\vec a_{\hat i}, b)$ holds for all $b\in \N$, which
    completes the proof of the claim.
    In particular, $R$ restricted to injective sequences
    is definable by a formula in the trivial language.

    Next, we need to deal with $R$ off of injective sequences.
    First note that if $\vec a$ is not injective, 
    then $\vec a$ is of the form $\vec \pi(\vec \rho(\vec a))$
    for some length $n-1$ sequence of projections $\vec \rho$ and some $(n-1)$-ary $\vec \pi$.  Indeed, if $a_i=a_j$, $i\neq j$,
    then $\vec \rho$ deletes the $a_j$ and
    $\vec \pi$ copies it back into the $j$th slot using $a_i$.
    Thus, off of injective sequences, $R$ is equivalent to
    \[
    \bigdoublevee_{\vec \rho} \ \ 
    \bigdoublevee_{\vec \pi} \Big [\vec a = \vec \pi(\vec\rho(\vec a)) \conj R_{\vec \pi}(\vec \rho(\vec a)) \Big],
    \]
    where $\vec \rho$ ranges over $n$-ary shuffling sequence
    of length $n-1$, and $\vec \pi$ ranges over $(n-1)$-ary
    shuffling sequences of length $n$.  
    So, it is enough to show each $R_{\vec \pi}$ above is definable
    with a formula in the trivial language.

    For a $(n-1)$-ary $\vec \pi$, $R_{\vec \pi}$ is an $(n-1)$-ary
    relation.  Each shuffle of $R_{\vec \pi}$ is also
    a shuffle of $R$, hence $R_{\vec \pi}$ also has the bad
    coding property.  By our inductive hypothesis,
    $R_{\vec \pi}$ is indeed definable with a formula
    in the trivial language.  
\end{proof}

For $x\in X_\mathcal{L}$, we will say that
$x$ has the \textbf{bad coding property} if 
$R_i^x$ has the bad
coding property for every $R_i\in \mathcal{L}$.

\begin{lemma}
\label{badcodinglemma2}
    Let $\mathcal{L}$ be a computable relational language
    and let $x\in X_\mathcal{L}$.  If $x$ has the bad coding
    property, then
    $x$ is the only structure in its $\cong_\mathcal{L}$-equivalence
    class.
\end{lemma}

\begin{proof}
    Since every $R_i^x$ has the bad coding property,
    they are all definable in the trivial language.
    It follows that every bijection is an automorphism of $x$.

    Suppose $y\cong_\mathcal{L} x$ via $f:\N\rightarrow \N$.
    Then, for every $R_i\in \mathcal{L}$ and any tuple $\vec a\in \N^{\text{arity}(R_i)}$,
    \[
    R_i^y(\vec a) \iff R_i^x(f(\vec a)) \iff R_i^x(\vec a).
    \]
    Thus, $x=y$.
\end{proof}

\begin{theorem}
\label{cisothm}
    Let $\mathcal{L}$ be a computable relational language.
    Then, computable isomorphism of $\mathcal{L}$-structures
    on $\N$
    is $\Pi^0_2$-graphable with diameter $2$.
\end{theorem}

\begin{proof}
    \ignore{Let $\text{BC}\subset X_\mathcal{L}$ be
    the set of structures which have the bad coding property.
    It is easy to compute that $\text{BC}$ is $\Pi^0_1$.}
    We begin by defining a predicate 
    $Q\subset X_\mathcal{L}\times \N^3$.
    When $Q(x, i, p, u)$ holds, $p$ codes a shuffling sequence
    $\vec \pi$ for $R_i^x$ of arity $k$, and $u$ codes a sequence $\vec a$
    of numbers whose length is $k-1$.\footnote{
    Note that if $k-1=0$, then $u$ is just the code of the empty sequence.
    }
    We will conflate the objects with the codes and just
    talk about $Q(x, i, \vec \pi, \vec a)$.  We will
    also use the computable wellordering 
    (of order type $\omega$) on these objects
    that comes from the numerical codes.
    $Q(x, i, \vec \pi, \vec a)$ holds when
    \begin{enumerate}[(Q1)]
    \item  $(i, \vec \pi)$ is least such that $(R_i^x)_{\vec{\pi}}$
    is a coding relation; and
    \item $\vec a$ is the least injective sequence so that there exist
    $b, c\in \N$, distinct from the elements of $\vec a$, 
    such that $(R_i^x)_{\vec{\pi}}(\vec a, b)$ and
    $\neg (R_i^x)_{\vec{\pi}}(\vec a, c)$ both hold.  
    \end{enumerate}
    A routine computation shows that $Q$ is $\Delta^0_2$.

    Let $f:\N\rightarrow \N$ be the total computable
    function defined in the proof of Theorem \ref{1equiv}.
    It will again provide us with bounds for searches
    for indices of computable isomorphisms.  

    Now, define $G\subset X_{\mathcal{L}}\times X_{\mathcal{L}}$
    as follows: $G(x, y)$ holds when
    \begin{enumerate}[(G1)]
    \item $y$ does not have the bad coding property; and
    \item $(\forall i, \vec \pi, \vec a, n)$ if
    $Q(x, i, \vec \pi, \vec a)$ holds and $n>0$ is the unique
    number in $\N\setminus \{\vec a\}$ such that $(R_i^x)_{\vec{\pi}}(\vec a, n)$ has a different truth
    value than all smaller numbers in $\N\setminus \{\vec a\}$, then ${(\exists e<f(n))}$ such that $\varphi_e$
    is an isomorphism from $x$ to $y$.
    \end{enumerate}
    It is easy to compute that $G$ is $\Pi^0_2$ and to see
    that $G(x, y)$ implies that $x$ and $y$ are computably isomorphic.
    What is left to show is that for any pair
    of distinct $x, x'\in X_\mathcal{L}$ with 
    $x\cong_\mathcal{L}^c x'$, there is a $y$ such that
    $G(x, y)$ and $G(x', y)$ both hold.

    Let $x, x'$ be distinct and computably isomorphic.
    Fix some $e\in \N$ so that $\varphi_e$ is an isomorphism
    from $x'$ to $x$.  
    Since they are distinct, it follows from Lemma 
    \ref{badcodinglemma2} that $x$ does not have the bad coding
    property.  Let $(i, \vec \pi)$ be least such that
    $(R_i^x)_{\vec{\pi}}$ is a coding relation.  Let
    $k+1$ be the arity of $\vec \pi$.\footnote{
    Note that if $k=0$, i.e., if $(R_i^x)_{\vec{\pi}}$ is unary,
    then the following construction still works, just ignoring
    the extra parameters.  In fact, the construction is
    exactly the same as the one in the proof of Theorem \ref{1equiv}.
    }
    Fix the least injective $k$-tuple 
    $\vec a=(a_0, \dots, a_{k-1})$ so that there exist
    $b, c\in \N$ such that $(R_i^x)_{\vec{\pi}}(\vec a, b)$ and
    $\neg (R_i^x)_{\vec{\pi}}(\vec a, c)$ both hold. 
    Fix $n$ which is larger than $k, e, b, c$, and 
    $\max(\vec a)+1$.  Define
    \begin{multline*}
    B:=\{m \in \N\setminus \{\vec a\} :m<2n \conj 
    (R_i^x)_{\vec{\pi}}(\vec a, m)\}, \\
    C:= \{m\in \N \setminus\{\vec a\}: 
    m<2n \conj \neg (R_i^x)_{\vec{\pi}}(\vec a, m)\}
    \end{multline*}
    By our choice of $n$, $B$ and $C$ are nonempty and
    disjoint, and
    $B\cup C$ has $2n-k$ many elements.  
    Clearly, at least one of them has $\geq n-k$ many elements.
    
    Now, 
    we will now describe how to build a new structure $y$
    which is the pushforward of $x$ by a finite swapping,
    using only transpositions of numbers
    $<2n$.  However, we will for now make two simplifying
    assumptions that we will discuss how to deal with later:
    \begin{enumerate}[(I)]
    \item we assume $B$ has $\geq n-k$ many elements; and
    \item we assume $\vec a=(0, 1, \dots, k-1)$.
    \end{enumerate}
    
    We describe all the needed transpositions to build $y$, under
    the simplifying assumptions (I) and (II).  Our main
    goal is to make $n$ satisfy (G2) in the structure $y$.  
    There are two cases. 

    \textit{Case 1}, $n\in C$.  Since $[k, n-1]\cap C$ has
    at most $n-k$ elements, by (I) we have enough elements of $B$
    to swap out all the elements of $[k, n-1]\cap C$ with elements of $B$.

    \textit{Case 2,} $n\in B$.  First, swap $n$ with some
    element of $[k, n-1]\cap C$ (say, $c$).  Now,
    $([k, n-1]\cap C)\setminus \{c\}$ has at most $n-k-1$
    many elements, so we can swap them all out with elements
    of $B\setminus \{n\}$, which has at least $n-k-1$ many
    elements.  

    The new structure $y$ has the following properties:
    \begin{enumerate}[(i)]
    \item $(i, \vec \pi)$ is least such that $(R^y_i)_{\vec \pi}$
    is a coding relation.  This is because isomorphisms do not
    change whether a definable relation in the structure
    is a coding relation.
    \item $\vec a= (0, 1, \dots, k-1)$ has the property that 
    there are $b, c\in \N$ such that $(R_i^y)_{\vec{\pi}}(\vec a, b)$
    and $\neg (R_i^y)_{\vec{\pi}}(\vec a, c)$ both hold. Moreover, as long
    as we have chosen a reasonable computable ordering of tuples,
    $\vec a= (0, 1, \dots, k-1)$ is the least such 
    injective $k$-tuple.
    \item \sloppy By construction, $n$ is the least number 
    different from
    $\vec a$ such that ${(R_i^x)_{\vec{\pi}}(\vec a, n-1)}$
    and 
    $(R_i^x)_{\vec{\pi}}(\vec a, n)$ have different
    truth values.
    \end{enumerate}
    Now, $x$ and $y$ are computably isomorphic by
    a finite swapping using transpositions of numbers $<2n$.
    Thus, $x$ and $y$ are isomorphic by a computable
    function with index $<f(n)$.  This implies $G(x, y)$ holds.
    Moreover, if you compose this function (on the left)
    with the isomorphism $\varphi_e$ from $x'$ to $x$,
    you have an isomorphism from $x'$ to $y$, again with
    index $<f(n)$.  Thus, $G(x', y)$ also holds.  

    What is left now is to describe how to do away with
    simplifying assumptions (I) and (II). If (I)
    fails, then we do the same construction, just switching
    the roles of $B$ and $C$.  If assumption (II) is not true,
    then we do our construction of $y$ in two stages.  In
    the first stage, we build $y'$ just by swapping each $a_i$
    with $i$.  All the $a_i$ were assumed to be less than $n$,
    so we have now a structure $y'$ on which we can apply our
    previous construction that we did under assumption (II).  
    This still results in only one
    step in the graph from $x$ to $y$ since we still only
    need finitely many transpositions of elements $<2n$.  
\end{proof}

We now explore several consequences of the theorem.
The following is immediate.

\begin{corollary}
Let $\mathcal{L}$ be a computable relational language
and let
$\mathcal{D}\subset X_\mathcal{L}$ be an isomorphism-invariant
class of $\mathcal{L}$-structures which is
$\Pi^0_2$.  Then, the restriction of
$\ciso$ to structures in $\mathcal{D}$,
\[
x (\cong_{\mathcal{L}}^{c}\restriction \mathcal{D}) y \iffdef 
x, y\in \mathcal{D} \conj x\ciso y
\]
is $\Pi^0_2$-graphable with diameter $2$.
In particular, computable isomorphism of linear orders
is $\Pi^0_2$-graphable with diameter $2$.
\end{corollary}

Just like we did for $m$-equivalence of sets,
we can also adapt the proof of Theorem \ref{cisothm}
to the case of computable biembeddability.
This equivalence relation on $X_\mathcal{L}$, denoted  
$\text{BE}_\mathcal{L}^c$, is defined so
that $x (\text{BE}_\mathcal{L}^c) y$ if and only
if there is a computable $\mathcal{L}$-embedding
from $x$ to $y$ and one from $y$ to $x$.  

\begin{corollary}
    For any computable relational language $\mathcal{L}$,
    the equivalence relation of computable biembeddability
    on $X_\mathcal{L}$
    is $\Pi^0_2$-graphable with diameter $2$.  
\end{corollary}

\begin{proof}
    The proof is nearly identical to the proof
    of Theorem \ref{cisothm}.  When defining the 
    graphing $G(x, y)$, just update (G2) so that $f(n)$
    provides an upper bound for $e_1$ and $e_2$
    which are indices for embeddings from $x$ to $y$
    and from $y$ to $x$, respectively.
\end{proof}

Next, we turn our attention to arbitrary computable languages,
which are of course allowed to have function symbols and constant
symbols in addition to relation symbols.  We will consider
constants symbols to just be symbols for nullary functions.
Just like with relational languages, a 
countable language $\mathcal{L}$ is \textbf{computable} if the assignment
of the symbols to their arity is computable.  The space of
$\mathcal{L}$-structures presented on $\N$ is then
\[
X_\mathcal{L}= \prod_i 2^{(\N^{\text{arity}(R_i)})} \times
\prod_i \N^{(\N^{\text{arity}(f_i)})},
\]
where $R_i$ are the relation symbols and $f_i$ are the functions symbols.
This space is a recursive Polish space, 
but, of course, it is no longer
compact when there are function symbols.  

\begin{corollary}
    Let $\mathcal{L}$ be a computable language.  Then,
    $\cong_\mathcal{L}^c$ is $\Pi^0_2$-graphable with
    diameter $2$.
\end{corollary}

\begin{proof}
    We will construct a relational language $\mathcal{L}'$
    such that $\cong_\mathcal{L}^c$ is invariantly computably
    reducible to $\cong_{\mathcal{L}'}^c$.  By Proposition 
    \ref{invariantprop} and Theorem \ref{cisothm}, this implies the
    result.  

    $\mathcal{L}'$ is obtained from $\mathcal{L}$ by replacing
    every function symbol $f\in \mathcal{L}$ with a a relation 
    symbol $R_f$ with arity $1+\text{arity}(f)$.  To build the invariant
    computable reduction $F:X_\mathcal{L} \rightarrow X_{\mathcal{L}'}$,
    it is enough to describe for every $x\in X_\mathcal{L}$ how 
    the structure
    $F(x)$ interprets all the symbols in $\mathcal{L}'$.  
    If $R$ is a relation symbol from $\mathcal{L}$ (and hence in 
    $\mathcal{L}')$, then the interpretation is unchanged, i.e.,
    $R^{F(x)}=R^x$. For a symbol $R_f$, where $f$ is a function symbol
    from $\mathcal{L}$, let $R_f^{F(x)}$ be the graph of $f^x$.  

    \sloppy $F$ is clearly computable and a reduction from $\cong_\mathcal{L}^c$ to
    $\cong_{\mathcal{L}'}^c$, since $\mathcal{L}$-isomorphisms
    preserve the graphs of the function symbols.
    To show it is invariant, it is enough to show that its image
    is closed under isomorphism.  
    Suppose $y\in X_{\mathcal{L}'}$ is isomorphic to
    $F(x)$ for some $x\in X_\mathcal{L}$.  Since being
    the graph of a function is a first-order property, it
    follow that each $R_f^y$ is
    the graph of a function.  Then, we can easily build an 
    $\mathcal{L}$-structure $y'$ with $F(y')=y$ by interpreting
    the relation symbols the same way they are interpreted in $y$
    and interpreting each function symbol $f$ in $y'$ so that
    its graph is $R_f^y$.  
\end{proof}

\section{Friedman-Stanley jumps}
\label{fsjumps}

If $E$ is an equivalence relation on $X$, then
its \textbf{Friedman-Stanley jump} is the 
equivalence relation $E^+$ on $X^\N$ defined by
\[
(x_i)E^+(y_i) \iff \{[x_i]_E : i\in \N\} = \{[y_i]_E : i\in \N\}.
\]
If $E$ is in $\Gamma$, then $E^+$ is in $\forall^\N\exists^\N \Gamma$.
In particular, if $E$ is arithmetical (respectively, Borel), then
$E^+$ is also arithmetical (resp., Borel).

For example, if $=_X$ denotes the equivalence relation
of equality on $X$, then $=_X^+$ is just
\[
(x_i)=_X^+ (y_i) \iff \{x_i:i\in\N\} = \{y_i : i\in \N\}.
\]
Both $=_{\N^\N}^+$ and $=_{\cantor}^+$ are $\Pi^0_3$,
while $=_\N^+$ is $\Pi^0_2$.  

The main theorem of this section gives a way to definably convert a
graphing of $E$ into a graphing of $E^+$.
The setting for this theorem will be a recursive Polish space
that has a $\Sigma^0_1$ definable strict linear ordering of
the space.  Such spaces include $\N$ and $\mathbb{R}$.  The
Kleene-Brouwer ordering gives a $\Sigma^0_1$ strict
linear ordering of $\cantor$, $\N^\N$, and $X_\mathcal{L}$
for any computable language $\mathcal{L}$.  Moreover, 
spaces which are computably isomorphic to products
of these spaces admit $\Sigma^0_1$ strict linear oderings.


\begin{theorem}
\label{fsjumpthm}
    Let $\Gamma$ be one of the pointclasses $\Sigma^0_n$, $n\geq 1$,
    or $\Pi^0_n$, $n\geq 2$.\footnote{
    Actually, the theorem also works for any pointclass $\Gamma$
    which contains $\Sigma^0_1$ and $\Pi^0_1$, 
    is closed under computable substitutions, $\vee$, $\&$, and computable substitutions, and 
    bounded existential quantification over $\N$.  In particular,
    the theorem holds for $\Delta^1_1$ and the Borel sets.
    }  
    Let $X$ be a recursive Polish space that has a 
    $\Sigma^0_1$ strict linear ordering $\prec$, and
    let $E$ be an equivalence relation on $X$.  
    If $E$ is $\Gamma$-graphable with finite diameter $\ell$,
    then $E^+$ is 
    $\forall^\N\Gamma$-graphable with diameter $\max(2, \ell)$.
\end{theorem}

\begin{proof}
    Fix a $\Gamma$ graphing $G$ of $E$ which has finite
    diameter $\ell$.  
    We introduce some notation.  Let $d_G(x, y)$ be the 
    distance between $x$ and $y$ in the graph $G$, i.e., 
    $d_G(x, y)$ is the 
    length of the shortest path in $G$ which connects them,
    if such a path exists, and $\infty$ otherwise.  
    For $x\in X$ and $(y_i)\in X^\N$, let
    $d_G(x, (y_i))=\min_{i\in\N} d_G(x, y_i)$.  
    Finally, let $d_G((x_i), (y_i))$ be the maximum
    among all the $d_G(x_j, (y_i))$ and $d_G(y_j, (x_i))$.
    We note a few basic properties:
    \begin{enumerate}[(i)]
    \item $d_G(x, y)\leq 1$ if and only
    if $x=y$ or $G(x, y)$.  In particular, the condition
    ``$d_G(x,y)\leq 1$'' is $\forall^\N \Gamma$.
    \item $(x_i)E^+(y_i)$ iff $d_G((x_i), (y_i))$ is finite iff 
    $d_G((x_i), (y_i))\leq \ell$,
    where $\ell$ is the diameter of $G$.
    \item If $(\tilde x_i)=_{X}^+ (x_i)$, then
    $d_G((\tilde x_i), (y_i)) = d_G((x_i), (y_i))$.
    \end{enumerate}
    
    Define the partial function 
    $g:X^\N\times \N\rightharpoonup \N$ by setting
    $g((x_i), n)\downarrow =i$ if and only if $i$ is the $(n+1)$st number
    such that $x_{2i+1}\prec x_{2i}$.  
    
    Set 
    $C:=\{(x_i): (\forall N)(\exists i>N)~x_{2i+1}\prec x_{2i}\}$,
    which is $\Pi^0_2$.  Since $\Gamma$ contains $\Sigma^0_1$,
    $C$ is in $\forall^\N \Gamma$.  
    If $(x_i)\in C$ then
    $n\mapsto g((x_i), n)$ is a total, strictly increasing function.
    Moreover, when $(x_i)\in C$, the condition $g((x_i), n)\downarrow = i$
    is $\Sigma^0_1$ in $(x_i)$.  
    
    There is a type of converse to the above claim.
    Given $(x_i)$ with $x_{2i}\neq x_{2i+1}$ for all $i$ 
    and any total,
    strictly increasing $f(n)$ then, 
    by swapping the appropriate $x_{2i}$
    with $x_{2i+1}$, we can form a rearrangement $(\tilde x_i)\in C$
    of $(x_i)$ such that $g((\tilde x_i), n) = f(n)$ for all $n\in\N$.

    Define $H((x_i), (z_i))$ so that it holds if
    either of the following hold:
    \begin{enumerate}[(H1)]
    \item $(\forall i)[d_G(x_i, z_i)\leq 1]$; or
    \item $(z_i)\in C$ and
    $(\forall n)(\exists j_0, j_1\leq g((z_i), n))[{d_G(x_n, z_{j_0}),
    d_G(x_{j_1}, z_n)\leq 1}]$.
    \end{enumerate}
    Using our assumptions about $\Gamma$, 
    it is routine to compute that $H$ is indeed $\forall^\N \Gamma$.
    Moreover, it is clear that $H((x_i), (z_i))$ implies
    $(x_i)E^+(z_i)$.  Note that in (H2), the condition $(z_i)\in C$
    ensures that the other conjunct cannot be vacuously true
    because of divergence of $g((z_i), n)$.  

    We isolate the main tools of the proof in
    the following two claims.

    \smallskip

    \textit{Claim 1.} If
    $d_G((x_i), (w_i)), d_G((y_i), (w_i)) \leq 1$, then
    there is a rearrangement $(z_i) \equiv_{X}^+ (w_i)$
    such that $H((x_i), (z_i))$ and $H((y_i), (z_i))$.
    In particular (when $(x_i)=(y_i)$), if $d_G((x_i), (w_i))\leq 1$,
    then there exists 
    $(z_i) \equiv_{\cantor}^+ (w_i)$
    such that $H((x_i), (z_i))$

    \smallskip

    \textit{Proof of Claim 1.} 
    If $w_i=w_j$ for all $i, j\in \N$, then 
    it is easy to see that $H((x_i), (w_i))$ and
    $H((y_i), (w_i))$ hold by virtue of (H1), hence we
    can take $(z_i):=(w_i)$.  If $(w_i)$ is not constant, than
    by replacing it with a $=_{X}^+$-equivalent sequence,
    we may assume $w_{2i}\neq w_{2i+1}$ for all $i\in \N$.
    Now, choose
    a strictly increasing $f(n)$ such that for every $n\in \N$
    there exists $j_0, j_1, j_0', j_1'<f(n)$ such that
    $d_G(x_n, w_{j_0})$, $d_G(x_{j_1}, w_n)$, $d_G(y_n, w_{j_0'})$,
    $d_G(y_{j_1'}, w_n)$ are all $\leq 1$.  We can rearrange $(w_i)$,
    by swapping pairs $w_{2i}, w_{2i+1}$ as needed, to get
    $(z_i)\in C$ such that $g((z_i),n) = f(n)$ 
    for all $n\in \N$.
    Note that every $z_i$ is equal to one of $w_{i-1}, w_i, w_{i+1}$;
    thus, if there is $j<f(n)$ so that $w_j$ has a certain property,
    then there is $j\leq f(n)$ so that $z_j$ has that same property.
    Using this, it is easy to see that both $H((x_i), (z_i))$
    and $H((y_i), (z_i))$ hold by virtue of clause (H2). 
    This completes the proof of Claim 1.

    \smallskip

    \textit{Claim 2.} \sloppy If $d_G((x_i), (y_i)) = k+1$, then
    there exists $(w_i)$ such that ${d_G((x_i), (w_i))=1}$ 
    and
    $d_G((y_i), (w_i)) = k$.

    \smallskip

    \textit{Proof of Claim 2.} For each $i\in \N$, if 
    $d_G(x_i, (y_j)) \leq k$, then
    just take $w_{2i}:= x_i$; otherwise, $d_G(x_i, (y_j))=k+1$
    and we can pick $w_{2i}$ so that $G(x_i, w_{2i})$ and
    $d_G(w_{2i}, (y_j))=k$.  Then, for $2i+1$, if
    $d(y_i, (x_j))\leq k$, take $w_{2i+1}:= x_j$ for some choice
    of $x_j$ with $d(y_i, x_j)\leq k$; otherwise, we can find
    some $x_j$ with $d(y_i, x_j)=k+1$ and 
    we pick $w_{2i+1}$ so that
    $G(x_j, w_{2i+1})$ and $d_G(w_{2i+1}, y_i)=k$.  It is
    easy to check that $d_G((x_i), (w_i)) = 1$ and that
    $d_G((y_i), (w_i)) =k$, finishing the proof of Claim 2. 

    \smallskip

    Note
    that if $d_G((x_i), (y_i)) \leq 1$, then we can
    apply Claim 1 with $(w_i):=(x_i)$, to get a path of length
    $2$ in $H$ connecting $(x_i)$ to $(y_i)$. 
    Next, we show by induction that for 
    $k\geq 2$, if $d_G((x_i), (y_i))=k$,
    then $(x_i)$ and $(y_i)$ are connected in $H$ by a
    path of length $k$.  
 
    If $d_G((x_i), (y_i))=2$, then we first apply Claim 2 to get
    $(w_i)$ with $d_G((x_i), (w_i))= d_G((y_i), (w_i)) = 1$.
    Then apply Claim 1 to get $(z_i)\equiv_{\cantor}^+ (w_i)$ 
    with $H((x_i), (z_i))$
    and $H((y_i), (z_i))$.  Finally, suppose 
    $d_G((x_i), (y_i))=k+1>2$.  We 
    apply Claim 2 to get $(w_i)$ with $d_G((x_i), (w_i))=1$
    and $d_G((y_i), (w_i)) = k$.  Then, by Claim 1, there
    is a $(z_i)=_{\cantor}^+ (w_i)$ 
    with $H((x_i), (z_i))$.  Note that we still have
    $d_G((y_i), (z_i))=k$ by (iii). By induction, there
    is a path of length $k$ in $H$ from $(z_i)$ to $(y_i)$.  
\end{proof}

The theorem has many consequences, the most immediate
of which applies to all the equivalence relations
we have already proved are $\Pi^0_2$-graphable
with diameter $2$.

\begin{corollary}
    All of the following equivalence relations have
    the property that their finite order Friedman-Stanley
    jumps are all $\Pi^0_2$-graphable with diameter $2$:
    $\equiv_T$, $\equiv_1$, $\equiv_m$, and
    $\cong_{\mathcal{L}}^c$,
    where $\mathcal{L}$ is a computable relational language.
    \ignore{and $\mathcal{C}\subset X_\mathcal{L}$ is a $\Pi^0_2$ isomorphism
    invariant class), and
    computable biembeddability
    of $\mathcal{L}$-structures on $\N$ (where $\mathcal{L}$ is
    a computable relational language).}  
\end{corollary}

Next, we point out a corollary that
comes from the fact that every equivalence relation $E$
is graphed by $G= {E\setminus (=_X)}$.  Moreover, this
graphing has diameter $\leq 1$.

\begin{corollary}
    Let $X$ be a recursive Polish space and with a 
    $\Sigma^0_1$ strict linear ordering, let $\Gamma$
    be a pointclass which contains $\Sigma^0_1$
    and is closed under computable substitutions, and
    let $E$ be an equivalence relation on $X$ which is in $\Gamma$. 
    Then, $E^+$ has a $\forall^\N \Gamma$-graphing
    of diameter $2$.  
\end{corollary}

Note this is a nontrivial statement since, in general, when $E$
is in $\Gamma$, $E^+$ is
in $\forall^\N\exists^\N \Gamma$.  In particular,
for every recursive Polish space $X$, the equality
equivalence relation $=_X$ is $\Pi^0_1$, hence
it is also $\Pi^0_2$.

\begin{corollary}
    If $X$ is a recursive Polish space with a $\Sigma^0_1$ strict
    linear ordering, then $=_X^+$ is $\Pi^0_2$-graphable with
    diameter $2$.  Moreover, all finite order
    Friedman-Stanley jumps of
    $=_X$ are also $\Pi^0_2$-graphable with diameter $2$.
\end{corollary}

The theorem also has consequences
for Borel graphability.
Note that every uncountable Polish space
is Borel isomorphic to $\mathbb{R}$, and so has a Borel linear ordering.

\begin{corollary}
    If $E$ is a Borel graphable with finite 
    diameter $\ell$, then $E^+$ is also
    Borel graphable with diameter $\max(2, \ell)$.
\end{corollary}

The above corollary requires that the
Borel graphing has finite diameter (which is quite
often the case).  We now show that one can drop this
requirement. Note, however, that the construction
in the following proof when applied
to a Borel graphing of $E$ of finite diameter $\geq 2$
will produce a Borel graphing of $E^+$ with 
larger diameter.

\begin{theorem}
    If $E$ is Borel graphable, $E^+$
    is also Borel graphable.
\end{theorem}

\begin{proof}
    Let $E^\N$ be the infinite product of $E$, i.e.,
    $E^\N$ is the equivalence relation on $X^\N$ defined
    by
    \[
    (x_i)E^\N(y_i) \iff (\forall i)~x_iEy_i.
    \]
    Note that $E^\N \subset E^+$.
    It is proved in \cite{akl} that Borel
    graphability is closed under infinite products.
    So, let $H$ be a Borel graphing of $E^\N$.

    Now, we define a graph $G$ on $X^\N$
    by letting $G((x_i), (y_i))$ hold when $(x_i)\neq (y_i)$ and
    either one of
    $(x_i)(=_X^+)(y_i)$ or $H((x_i), (y_i))$ holds. 
    It is clear that $G$ is Borel and $G((x_i), (y_i))$
    implies $(x_i)E^+(y_i)$.  

    Suppose $(x_i)$ and $(y_i)$ are distinct and $E^+$-equivalent.  Now, define $(\tilde x_i)$ and $(\tilde y_i)$
    as follows.  For each $n\in \N$, choose $i_n$ and $j_n$
    so that $x_nEy_{i_n}$ and $x_{j_n}Ey_n$.  
    Now, set $\tilde x_{2n} := x_n$, 
    $\tilde x_{2n+1}:= x_{j_n}$, $\tilde y_{2n} := y_{i_n}$
    and $y_{2n+1} := y_n$.  Immediately, we have
    $(x_i)=_X^+ (\tilde x_i)$, $(y_i)=_X^+ (\tilde y_i)$
    (so that both pairs are $G$-adjacent), and
    that $(\tilde x_i) E^\N (\tilde y_i)$.  So,
    we can use an $H$-path from $(\tilde x_i)$ to
    $(\tilde y_i)$ to create a $G$-path from $(x_i)$
    to $(y_i)$. 
\end{proof}

\ignore{
\section{Questions}

All of our graphing results gave us graphings of diameter $2$,
which motivates the following question.

\begin{question}
    Is there an arithmetical equivalence relation $E$
    and a pointclass $\Gamma$ so that $E$ has a $\Gamma$-graphing
    $G$ of diameter $\geq 3$ but no $\Gamma$-graphing
    of diameter $2$?  Similarly, is there an arithmetical $E$
    which has $\Gamma$-graphing of infinite diameter but
    no $\Gamma$-graphing of finite diameter?
\end{question}

Lower bound results on the definability of graphings seem difficult,
but would be necessary to prove that the graphings obtained
here are optimal.  

\begin{question}
    Is $\equiv_1$, as an equivalence relation on $\cantor$,
    $\Sigma^0_2$-graphable?
\end{question}
}

\bibliographystyle{alpha}
\bibliography{main.bib}

\end{document}